\documentclass[12pt]{article}
\usepackage{amsmath}
\usepackage{amssymb}
\usepackage{amsfonts, amsthm}
\begin{document}
\title{An Eulerian-Lagrangian Approach to the Navier-Stokes
equations.}
\author{Peter Constantin
\\Department  of Mathematics\\The University of Chicago }
\maketitle
\newtheorem{thm}{Theorem}
\newtheorem{prop}{Proposition}

\section{Introduction}

This work presents an Eulerian-Lagrangian approach to the Navier-Stokes  
equation. An Eulerian-Lagrangian
description of the Euler equations has been
used in (\cite{cel}, \cite{celw}) for local existence results and
constraints on blow-up. Eulerian coordinates 
(fixed Euclidean coordinates) are natural for 
both analysis and laboratory experiment. Lagrangian variables
have a certain theoretical appeal. In this work I present an approach 
to the Navier-Stokes equations that is phrased in
unbiased Eulerian coordinates, yet describes objects that have  
Lagrangian significance: particle paths, their dispersion
and diffusion. The commutator between Lagrangian and Eulerian 
derivatives plays an important role in the Navier-Stokes equations:  
it contributes a singular perturbation to the Euler equations, in addition
to the Laplacian. The Navier-Stokes equations are shown to be equivalent to
the system
$$
\Gamma v = 2\nu C\nabla v
$$
where $C$ are the coefficients of the commutator between Eulerian and 
Lagrangian derivatives, and  $\Gamma$ is the operator of material derivative 
and viscous diffusion. The physical pressure is not explicitly present in this
formulation. The Eulerian velocity $u$ is related to $v$ in a non-local 
fashion,
and one may recover the physical pressure dynamically from the evolution of
the gradient part of $v$. When one sets $\nu =0$ the commutator
coefficients $C$ do not enter the equation, and then $v$ is a passive
rearrangement of its initial value (\cite{serr}, \cite{gold}, \cite{goldet},
\cite{hunt}). When $\nu\neq 0$ the perturbation  involves the curvature of 
the particle paths, and the gradients of
$v$: a singular perturbation. Fortunately, the coefficients $C$ start from 
zero, and, as long as they remain small $v$ does not grow too much.

A different but not unrelated approach (\cite{kuz}, (\cite{ose}) is based on a 
variable $w$ that has the same curl as the Eulerian velocity $u$. The
velocity is recovered then from $w$ by projection on divergence-free
functions. The evolution equation for $w$ 
$$
\Gamma w + \left(\nabla u\right )^*w = 0,
$$ 
conserves local helicity and circulation (when $\nu=0$). We will refer to
this equation informally as ``the cotangent equation'' because it is the 
equation obeyed by the Eulerian gradient of any scalar $\phi$ that 
solves $\Gamma\phi = 0$. The variable $w$ is related to $v$: 
$$
w = (\nabla A)^*v
$$ 
where $A$ is the ``back-to-labels'' maps that corresponds when $\nu =0$
to the inverse of the Lagrangian path map. $A(x,t)$ is an active vector 
obeying 
$$
\Gamma A = 0,
$$
$A(x,0) =x$. Both $v$ and $A$ have a Lagrangian meaning when $\nu=0$, but
the dynamical development of
$w$ is the product of two processes, the growth
of the deformation tensor (given by the evolution of $\nabla A$) and the 
rearrangement of a fixed function, given by the evolution of $v$. 
In the presence of viscosity, $v$'s evolution is not by rearrangement only.
It is therefore useful to study separately the 
growth of $\nabla A$ and the shift of $v$.

Recently certain model equations (alpha-models) have been 
proposed (\cite{h}, \cite{mar})
as modifications of the Euler and Navier-Stokes equations. They
can be obtained in the context described above
simply by smoothing $u$ in the cotangent equation. 
Smoothing means that one replaces the linear zero-order 
nonlocal operator $u = {\mathbf P}w$ that relates $u$ to $w$ in the cotangent 
equation  by a smoothing operator, $u =K_{\alpha}{\mathbf P}w$.
When $\nu =0$ the models have a Kelvin circulation theorem. 
They cannot be models for Eulerian averaged equations, the equations 
describing mean flow in turbulence theory. These Reynolds 
equations are not conservative in the limit of zero viscosity: the 
fluctuations introduce additional stresses, the so-called Reynolds stresses, 
that preclude conservation of circulation along the mean flow. 
This is a simple yet fundamental objection to the identification
of alpha-models with Reynolds equations. The alpha-models might 
be models of Lagrangian averaged equations, an entirely
different concept.  (Operationally one may think of 
Eulerian averaging as long time averages at fixed position, and of 
Lagrangian averaging as long time averages at fixed initial label. However, 
a fixed initial label has no obvious physical meaning when one deals with 
ensembles of flows.)   

In this paper we consider the Navier-Stokes equations and
obtain rigorous bounds for the particle paths and for the virtual 
velocity $v$. The main bounds concern the Lagrangian displacement, its first
and second spatial derivatives, are obtained under general 
conditions and require no assumptions. Higher derivatives
can be bounded also under certain natural quantitative
smoothness assumptions.

Some of our bounds can be interpreted as a connection to
the Richardson pair dispersion law, one of the
empirical laws of fully developed turbulence (\cite{moya}), that is
consistent with the Kolmogorov two-thirds law (\cite{frisch}). The pair 
dispersion law states that the separation
$\delta $ between fluid particles obeys
$$
\left < |\delta |^2\right > \sim \epsilon t^3
$$
where $\epsilon$ is the rate of dissipation of energy and $t$ is time. This
is supposed to hold in an inertial range, for times $t$ that are neither too 
small (when the separation is ballistic) nor too large, 
when viscous and boundary effects are important. The law can be guessed by 
dimensional analysis by requiring the answer to depend solely on time and 
$\epsilon$. Precise 
laboratory Lagrangian experiments have recently begun to be capable of 
addressing Lagrangian quantities
with preliminary results that seem to be consistent with the 
Richardson law in some ranges (\cite{labexp}). 
If one considers the problem of estimating the pair dispersion mathematically
one is faced with the difficulty that the prediction seems to require both
non-Lipschitz, H\"{o}lder continuous velocities and Lagrangian particle paths. 
Our approach allows a rigorous
formulation and an upper bound 
$$
\left < |\delta |^2\right > \le 3|\delta_0 |^2 + 24E_0 t^2 + 6\epsilon_B t^3,
$$
that includes a reference to an initial displacement $\delta_0$, an initial
kinetic energy $E_0$ and a rigorous upper bound on $\epsilon$. The prefactors 
are probably not optimal. The conditions under which such a bound can be obtained are 
quite general, and there are no assumptions. In this paper we chose the case of periodic boundary conditions,
and body forces that have a characteristic length scale that remains finite
as the size of the periodic box is allowed to diverge. The bound 
$\epsilon_B$ does not depend on the size of box.
In many physically realistic situations one 
injects energy at the boundary; in that case
one can find $\epsilon_B$ independently of viscosity (\cite{dc}), without
any assumptions.

\section{Velocity and displacement}
The Eulerian velocity 
$u(x,t)$ has three components $u^i, \,\, i=1,\,2\,,  3$
and is a function of three Eulerian space coordinates
$x$ and time $t$. We 
decompose the Eulerian velocity 
$u(x,t)$:
\begin{equation}
u^i(x,t) = \frac{\partial A^m(x,t)}{\partial x_i}v_m(x,t) - 
 \frac{\partial n(x,t)}{\partial x_i}.\label{u}
\end{equation}
Repeated indices are summed. 
There are three objects that appear in this formula. The first one, $A(x,t)$,
has a Lagrangian interpretation. In the absence of viscosity, $A$ is the ``back-to-labels'' map, the inverse of the 
particle trajectory map $a\mapsto x= X(a,t)$. 
The vector
\begin{equation}
\ell(x,t) = A(x,t) - x\label{d}
\end{equation}
will be called the ``Eulerian-Lagrangian displacement vector'', or 
simply ``displacement''. $\ell$ joins the current Eulerian position $x$ to 
the original Lagrangian  position $a =A(x,t)$. $A(x,t)$ and $\ell (x,t)$
have dimensions of length, $\nabla A$ is non-dimensional. A pair of points, $a = A(x,t)$, $b = A(y,t)$
situated at time $t=0$ at distance $\delta_0 = |a-b|$ become separated by
$\delta_t = |x-y|$ at time $t$. From the triangle inequality it follows
that 
\begin{equation}
(\delta_t)^2  \le 3|\ell(x,t)|^2 + 3|\ell (y,t)|^2 + 3(\delta_0)^2.\label{tri}
\end{equation}
The displacement
can be used in this manner to bound pair dispersion. 

The second object in (\ref{u}), $v(x,t)$, has dimensions of velocity and, 
in the absence of viscosity, is just the initial velocity 
composed with the back-to-labels map (\cite{serr}, \cite{gold}, 
\cite{goldet}, \cite{hunt}). 
We call $v$ the ``virtual velocity''. Its evolution  marks the difference
between the Euler and Navier-Stokes equations most clearly.
The third object 
in (\ref{u}) is a scalar function $n(x,t)$ 
that will be referred to as ``the Eulerian-Lagrangian 
potential''. It plays a mathematical
role akin to that played by the physical pressure but has dimensions
of length squared per time, like the kinematic viscosity. 
If $A(x,t)$ is known, then there are four functions entering the decomposition
of $u$, three $v$-s and one $n$. If the velocity is divergence-free
$$
\nabla\cdot u = 0,
$$
then there is one relationship between the four unknown functions.

\section{Eulerian-Lagrangian derivatives and commutators}
When one considers the map $x\mapsto A(x,t)$ as a change of variables one
can pull back the Lagrangian differentiation with respect to particle
position and write it in Eulerian coordinates using the 
chain rule. Let us call this pull-back of Lagrangian derivatives
the Eulerian-Lagrangian derivative,
\begin{equation}
\nabla_A = Q^*\nabla_E.\label{la}
\end{equation}
Here 
\begin{equation}
Q(x,t) = \left (\nabla A(x,t)\right )^{-1},\label{Q}
\end{equation}
and the notation $Q^*$ refers to the transpose of the matrix $Q$. The
expression of $\nabla_A$ on components is
\begin{equation}
\nabla_A^i = Q_{ji}\partial_j\label{lai}
\end{equation} 
where we wrote  
$\partial_i$ for differentiation 
in the $i$-th Eulerian Cartesian coordinate direction,
$$
\partial_i = \nabla_E^i.
$$ 
The Eulerian spatial derivatives can be expressed in terms of the 
Eulerian-Lagrangian
derivatives via
\begin{equation}
\nabla_E^i = \left (\partial_i A_m\right )\nabla_A^m
\label{ial}
\end{equation}
The commutation relations
$$
[\nabla_E^i, \nabla_E^k] = 0, \quad \left [\nabla_A^i, \nabla_A^k\right ] = 0
$$
hold. 
The commutators between Eulerian-Lagrangian and Eulerian derivatives do
not vanish, in general:
\begin{equation}
\left [\nabla_A^i, \nabla_E^k\right ] = C_{m,k; i}\nabla_A^m.\label{str}
\end{equation}
The coefficients $C_{m,k;i}$ are given by 
\begin{equation}
C_{m,k;i} = \{\nabla_A^{i}(\partial_k \ell_m)\}.\label{cmki}
\end{equation}
Note that
$$
C_{m,k;i} = Q_{ji}\partial_i\partial_k A_m =
\nabla_A^i(\nabla_E^k A_m) = [\nabla_A^i,\nabla_E^k]A_m.
$$
These commutator coefficients play an important role in dynamics.

\section{The evolution of A}
We associate to a given divergence-free
velocity $u(x,t)$ the 
operator
\begin{equation}
\partial_t + u \cdot\nabla -\nu\Delta = {\Gamma_{\nu} (u,\nabla )}\label{L}.
\end{equation}
We write $\partial_t$ for time derivative.
We write $\Gamma$ for $\Gamma_{\nu}(u,\nabla)$ when the $u$ 
we use is clear from the context. The coefficient $\nu >0$ is the 
kinematic viscosity of the fluid. When applied to a vector or a matrix, $\Gamma$ acts
as a diagonal operator, i.e. on each component separately. The operator
$\Gamma$ obeys a maximum principle: 
If a function $q$ solves
$$
\Gamma  q = S
$$
and the function $q$ has homogeneous
Dirichlet or periodic boundary conditions, then the sup-norm $\|q\|_{L^{\infty}(dx)}$ satisfies
$$
\|q(\cdot, t)\|_{L^{\infty}(dx)} \le \|q(\cdot ,t_0)\|_{L^{\infty}(dx)}
+ \int_{t_0}^t\|S(\cdot, s)\|_{L^{\infty}(dx)}ds
$$
for any $t_0\le t$. The operator ${\Gamma_{\nu} (u,\nabla)}$ is not a derivation 
(that means an operator that satisfies the product rule); 
$\Gamma$ satisfies a product rule that is similar to that of a 
derivation: 
\begin{equation}
{\Gamma}(fg) = ({\Gamma }f)g + f({\Gamma}g) - 2\nu(\partial_k f)(\partial_k g)\label{lfg}.
\end{equation}

We require the 
back-to-labels map $A$ to obey
\begin{equation}
\Gamma A = 0. \label{aeq}
\end{equation}
By (\ref{aeq}) we express therefore the advection  {\em and diffusion}
of $A$. We will use sometimes the equation obeyed by $\ell $ 
\begin{equation}
\left (\partial_t   + u\cdot\nabla  -\nu \Delta \right )\ell 
+ u = 0
\label{deltaeq}
\end{equation}
which is obviously equivalent to (\ref{aeq}).
We will discuss  periodic boundary conditions
$$
\ell(x+Le_j,t) = \ell (x,t),
$$
where $e_j$ is the unit vector in the $j$-th direction. 
Some of our inequalities will  hold 
also for the physical boundary condition that require
$\ell (x,t) = 0$ at the boundary.

It is important to note that the initial data for the displacement is zero:
\begin{equation}
\ell (x,0) = 0.\label{deltain}
\end{equation}
The matrix $\nabla A(x,t)$ is invertible as long as the
evolution is smooth. This is obvious when $\nu =0$ because the determinant
of this matrix equals $1$ for all time, but in the viscous case 
the statement needs proof. We differentiate 
(\ref{aeq}) in order to obtain the equation obeyed by $\nabla A$
\begin{equation}
{\Gamma}(\nabla A)  + (\nabla A)(\nabla u)
= 0.\label{nablaeq}
\end{equation}
The product $ (\nabla A)(\nabla u)$ is matrix product in the order indicated. We consider

\begin{equation}
{\Gamma }Q  = (\nabla u)Q +  
2\nu Q\partial_k(\nabla A)\partial_k Q.\label{Qeq}
\end{equation}
It is clear that the solutions of both (\ref{nablaeq}) and
(\ref{Qeq}) are smooth as long as the advecting velocity $u$ is 
sufficiently smooth. It is easy to verify using (\ref{lfg}) that the matrix $Z=(\nabla A)Q -I$ obeys the equation
$$
\Gamma Z = 2\nu Z\partial_k(\nabla A)\partial_k Q
$$
with initial datum $Z(x,0) = 0$. Thus, as long as $u$ is smooth, $Z(x,t) = 0$ and it follows
that the solution $Q$ of (\ref{Qeq}) is the 
inverse of $\nabla A$.

The commutator coefficients $C_{m,k;i}$ enter the important
commutation relation between
the Eulerian-Lagrangian label derivative and $\Gamma$:
\begin{equation}
\left [{\Gamma}, \nabla_A^{i}\right ] = 2\nu C_{m,k ;i}
\nabla_E^k \nabla_A^{m}
 \label{lagcom}
\end{equation}
The proof of this formula can be found in Appendix B.

The evolution of the coefficients $C_{m,k;i}$
defined in (\ref{cmki}) can be computed 
using (\ref{nablaeq}) and (\ref{lagcom}):
$$
{\Gamma}\left (C_{m,k;i} \right ) =
 - (\partial_lA_m)\nabla_A^i(\partial_k(u_l))
$$
\begin{equation}
-(\partial_k(u_l))C_{m,l;i}+ 2\nu C_{j,l;i}\cdot\partial_l\left (C_{m,k;j}\right ).\label{curveq}
\end{equation}
The calculation leading to (\ref{curveq}) is presented in Appendix B.

\section{The evolution of v}
We require the virtual velocity to obey 
\begin{equation}
\Gamma v = 2\nu C\nabla v + Q^*f.\label{psieq}
\end{equation}
This equation is, on components
\begin{equation}
\Gamma_{\nu}(u,\nabla)v_i = 2\nu C_{m,k;i}\partial_kv_m + Q_{ji}f_j.\label{veq}
\end{equation}
The vector
$f=f(x,t)$ represents the body forces.
The boundary conditions are periodic
$$
v(x+Le_j,t) = v(x,t)
$$
and the initial data are, for instance
\begin{equation}
v(x,0) = u_0(x).
\label{psin}
\end{equation}
The reason for requiring the equation (\ref{psieq}) is
\begin{prop}
Assume that $u$ is given by the expression
 (\ref{u}) above and that the displacement
$\ell$ and the virtual velocity $v$ obey the equations (\ref{deltaeq}) and respectively (\ref{psieq}). Then the velocity
$u$ satisfies the Navier-Stokes equation
$$
\partial_t u + u\cdot\nabla u - \nu\Delta u + \nabla p = f
$$
with pressure $p$ determined from the Eulerian-Lagrangian potential by
$$
{\Gamma_{\nu} (u,\nabla)} n + \frac{|u|^2}{2}  + c = p 
$$
where $c$ is a free constant.
\end{prop}

\noindent{\bf Proof.} We denote for convenience
\begin{equation}
D_t = \partial_t + u\cdot\nabla .\label{dt}
\end{equation}
We apply $D_t$ to the velocity representation (\ref{u})
and use the commutation relation
\begin{equation}
\left [D_t, \partial_k\right ]g = -(\nabla u)^*\nabla g.\label{commut}
\end{equation}
We obtain
$$
D_t(u^i) = \left (\partial _i (D_tA^m)\right )v_m + (\partial_iA^m)
D_tv_m - \partial_i\left (\frac{|u|^2}{2} + D_t n\right ).
$$
We substitute the equations for $A$ (\ref{deltaeq}) and for $v$
(\ref{psieq}):
$$
D_t(u^i) = - \partial_i\left (\frac{|u|^2}{2} + D_t n\right ) +
 \left (\partial _i (\nu\Delta A^m)\right )v_m + 
$$
$$
(\partial_iA^m)\left \{\nu \Delta v_m +  Q^*_{mj}\left (2\nu\partial_k(\nabla\ell)^*_{jl}\partial_k v_l +f_j\right )\right\}.
$$
Now we use the facts that
$$ 
(\partial_iA^m)Q^*_{mj} =\delta_{ij}
$$
(Kronecker's delta), and 
$$
\partial_k(\nabla\ell )^*_{il} = \partial_k(\nabla A)^*_{il}
= \partial_k(\partial_i A^l)
$$
to deduce
$$
D_t(u^i) =  - \partial_i\left (\frac{|u|^2}{2} + D_t n\right ) + f_i 
$$
$$
+ \nu (\Delta \partial_iA^m )v_m + \nu (\partial_iA^m)\Delta v_m
+ 2\nu \partial_k(\partial_i A^l)\partial_k v_l 
$$
and so, changing the dummy summation index $l$ to $m$ in the last expression
$$
D_t(u^i) = - \partial_i\left (\frac{|u|^2}{2} + D_t n\right ) 
 + \nu\Delta ((\partial_i A^m)v_m) + f_i.
$$
Using (\ref{u}) we obtained
$$
D_t(u^i) =
\nu\Delta u_i - \partial_i\left (\frac{|u|^2}{2} -\nu\Delta n + D_t n\right ) + f_i
$$
and that concludes the proof.

\vspace{.5cm}

\noindent{\bf Observation} The incompressibility of velocity has not yet been 
used. This is why no restriction on
the potential 
$n(x,t)$ was needed. The incompressibility 
\begin{equation}
\nabla\cdot u = 0\label{div}
\end{equation}
can be imposed in two ways. The first approach is static: one considers the
ansatz (\ref{u}) and one requires that $n$ maintains the
incompressibility at each instance of time. This results in the equation
\begin{equation}
\Delta n = \nabla \cdot\left (\nabla A)^*v\right ).\label{incs}
\end{equation}
In this way $n$ is computed from $A$ in a time independent manner and
the basic formula (\ref{u}) can be understood as
\begin{equation}
u = {\mathbf P}\left ((\nabla A)^*v\right )\label{pu}
\end{equation}
where ${\mathbf P}$ is the Leray-Hodge projector on divergence-free functions.
The second approach is dynamic: one computes the physical 
Navier-Stokes
pressure 
\begin{equation}
p = R_iR_j(u^iu^j) + c\label{p}
\end{equation}
where $c$ is a free constant and $R_i = (-\Delta )^{-\frac{1}{2}}\partial_i$ is the Riesz transform for periodic boundary conditions.
The formula for $p$ follows by taking the divergence of the 
Navier-Stokes equation and using (\ref{div}). 
Substituting (\ref{p}) in the expression for the pressure in Proposition 1
one obtains the evolution equation
\begin{equation}
{\Gamma }n  = R_iR_j(u^iu^j) - \frac{|u|^2}{2} + c
\label{incd}
\end{equation}
for $n$. Incompressibility can be enforced either by solving at
each time the static equation (\ref{incs}) or by evolving $n$ according
to (\ref{incd}).
\begin{prop}
Let  $u$ be given by (\ref{u}) and assume that the displacement 
solves (\ref{deltaeq}) and that the virtual velocity solves
(\ref{psieq}). Assume in addition that the potential 
obeys (\ref{incs}) (respectively (\ref{incd})). Then $u$ obeys the incompressible Navier-Stokes equations, 
$$
\partial_t u + u\cdot\nabla u - \nu\Delta u + \nabla p = f, \quad \nabla\cdot u =0,
$$ 
the pressure $p$ satisfies (\ref{p}) and the potential obeys
also (\ref{incd}) (respectively (\ref{incs})).
\end{prop} 

The same results hold for the case of the whole ${\mathbf R}^3$
with boundary conditions requiring $u$ and $\ell$ to vanish at infinity.
In the presence of boundaries, if the boundary conditions for $u$ are homogeneous Dirichlet
($u=0$) then the boundary conditions for $v$ are 
Dirichlet, but not homogeneous. In that case one needs to solve either
one of the equations (\ref{incs}),(\ref{incd}) for $n$ 
(with Dirichlet or other physical boundary
condition)  and the $v$ equation (\ref{psieq}) with
$$
v = \nabla_A n
$$
at the boundary.

\begin{prop}
Let $u$ be an arbitrary spatially periodic
smooth function and assume that a displacement
$\ell$ solves the equation (\ref{deltaeq}) and a virtual velocity $v$ obeys 
the equation (\ref{psieq}) with periodic boundary conditions and 
with $C$ computed using $A = x+\ell$. Then $w$ defined by 
\begin{equation}
w_i = (\partial_iA^m)v_m
\label{w}
\end{equation}
obeys the cotangent equation
\begin{equation}
{\Gamma }w + (\nabla u)^*w = f.\label{weq}
\end{equation}
\end{prop}

\noindent{\bf Proof.} The proof is a straightforward calculation. One uses
(\ref{lfg}) to write
$$
\Gamma w_ i = (\partial_iA^m)\Gamma v_m + v_m\Gamma(\partial_iA^m) -
2\nu (\partial_k\partial_iA^m)\partial_kv_m.
$$
The equation (\ref{veq}) is used for the first term and the
equation (\ref{nablaeq}) for the second term. One obtains
$$
{\Gamma }w_i = f_i - (\partial_i u_j)w_j +2\nu\left \{(\partial_iA^m)C_{r,q; m}
\partial_q v_r -  (\partial_k\partial_iA^m)\partial_kv_m\right\}.
$$
The proof ends by showing that the term in braces vanishes because of 
the identity
$$
(\partial_iA^m)C_{r,q; m} = \partial_q\partial_iA^r.
$$

\vspace{.5cm}

An approach to the Euler equations based entirely on a variable
$w$ (\cite{kuz}, (\cite{ose}) is well-known. 
The function $w$ has the same curl as $u$, $\omega = \nabla\times u = \nabla\times w$. In the case of zero viscosity and no forcing, the local 
helicity $w\cdot\omega$ is conserved $D_t(w\cdot\omega) = 0$; this is 
easily checked using the fact that the
vorticity obeys the ``tangent'' equation $D_t\omega = (\nabla u)\omega$ 
and the inviscid, unforced  form of (\ref{weq}). The same proof verifies
the Kelvin circulation theorem
$$
\frac{d}{dt}\oint_{\gamma(t)} w\cdot dX = 0
$$ 
on loops $\gamma(t)$ advected by the flow of $u$. 
Although obviously related, the two variables $v$ and $w$
have very different analytical merits. While the growth of $w$
is difficult to control, in the inviscid case $v$ does not grow at all, and 
in the viscous case its growth is determined  by the magnitude
of $C$ which starts from zero. This is why we emphasize $v$ as
the primary variable and consider $w$ a derived variable.

\section{Gauge Invariance}
Consider a scalar function $\phi$. 
If one transforms $v\mapsto \tilde{v}= v + \nabla_A\phi$ and $n\mapsto \tilde{n} =
n+\phi$ then $u$ remains unchanged in (\ref{u}): $u\mapsto u$.
The requirement that $\nabla_E\cdot u = 0$ does not specify this arbitrary 
$\phi$. 

Assume now that the scalar $\phi$ is advected passively by $u$ and diffuses 
with diffusivity $\nu$:
$$
\Gamma \phi  = 0.
$$
Then, in view of (\ref{lagcom}), if $v$ solves (\ref{psieq}) then 
$$
\tilde{v} = v + \nabla_A \phi
$$
also solves (\ref{psieq}). If $n$ solves (\ref{incd}) then
$$
\tilde{n} = n + \phi
$$
also solves (\ref{incd}). If $w$ solves the equation
(\ref{weq}) then
$$
\tilde{w} = w + \nabla_E \phi
$$
also solves (\ref{weq}).  
The vector fields obtained
by taking the Eulerian gradient of passive scalars are homogeneous
solutions of (\ref{weq}). The vector fields obtained by taking
the Eulerian-Lagrangian gradient of passive scalars, $\nabla_A \phi$ 
are homogeneous solutions of (\ref{psieq}).  This can be used to show
that if one chooses an initial datum for $v$ that differs 
from $u_0$ by the gradient of an arbitrary function $\phi_0$ there is no 
change in the evolution of $u$.
\begin{prop}
Let each of two functions $v_j$, $j=1,2$ solve the system
$$
\Gamma (u_j,\nabla )v_j= 2\nu C_j\nabla v_j + Q_j^*f
$$
with periodic boundary conditions, coupled with
$$
\Gamma (u_j,\nabla )A_j = 0
$$
with periodic boundary conditions for $\ell_j = A_j -x$. Assume that the initial
data for $A_j$ are the same, $\ell_j (x,0) = 0$. Assume that each velocity
is determined from its corresponding virtual velocity by the rule
$$
u_j = {\mathbf P}\left ((\nabla A_j)^*v_j\right).
$$
Assume, moreover, that at time \, $t=0$ the virtual velocities differ by a
gradient
$$
{\mathbf P}v_1 ={\mathbf{P}}v_2 = u_0.
$$
Then, as long as one of the solutions $v_j$ is smooth one has
$$
u_1(x,t) = u_2(x,t),\quad A_1(x,t) = A_2(x,t)
$$
\end{prop}
The same kind of result can be proved for (\ref{weq}) using the Eulerian
gauge invariance.
\section{K-bounds}
We are going to describe here bounds that are based solely on the
kinetic energy balance in the Navier-Stokes equation ( (\cite{cfbook}) and references therein). These are
very important, as they are the only unconditional bounds that
are known for arbitrary time intervals. We call them kinetic energy bounds or in
short, K-bounds.  We start with the most important, the energy balance
itself. From the Navier-Stokes equation one obtains the bound
\begin{equation}
\int |u(x,t)|^2dx + \nu\int_{t_0}^{t}\int|\nabla u(x,s)|^2dx ds
\le K_0\label{en}
\end{equation}
with
\begin{equation}
K_0 = \min{\{ k_0; k_1\}}\label{K0}
\end{equation}
where
\begin{equation}
k_0 = 2\int|u(x,t_0)|^2dx + 3(t-t_0)\int_{t_0}^{t}
\int |f(x,s)|^2dx ds\label{k0}
\end{equation}
and
\begin{equation}
k_1 = \int|u(x,t_0)|^2dx + \frac{1}{\nu} \int_{t_0}^t\int |\Delta^{-\frac{1}{2}} f(x,s)|^2dxds
\label{k1}
\end{equation}
Note that we have not normalized the volume of the domain. The prefactors 
are not optimal. The energy balance holds for all
solutions of the Navier-Stokes equations. 
We took an arbitrary starting time $t_0$. The bound $K_0$ is a nondecreasing
function of $t-t_0$. We will use this fact tacitly below. 
In order to give a physical interpretation to this general bound it is useful to denote by
$$
\epsilon (s) = \nu L^{-3}\int |\nabla u (x,s)|^2dx
$$
the volume average of the instantaneous energy dissipation rate,
by
$$
E(t) =  \frac{1}{2L^3}\int|u(x,t)|^2dx
$$
the volume average of the kinetic energy; for any time dependent function $g(s)$, we write
$$
\left < g (\cdot )\right>_t = \frac{1}{t-t_0}\int_{t_0}^tg(s)ds
$$
for the time average. We also write
$$
F^2 = \left <L^{-3}\int |f(x,\cdot)|^2dx\right>_t,
$$
$$
G^2  = \left <L^{-3}\int |\Delta^{-\frac{1}{2}}f(x,\cdot)|^2dx\right>_t 
$$
and define the forcing length scale by
$$
L_f^2 = \frac{G^2}{F^2}.
$$
Then
(\ref{en}) implies
\begin{equation}
2 E(t) + (t-t_0)\left < \epsilon(\cdot )\right >_t \le
4E(t_0) + (t-t_0)F^2\min{\left \{\frac{L_f^2}{\nu}; 3(t-t_0)\right \}}. 
\label{epsbound}
\end{equation}
After a long enough time
$$
t -t_0 \ge \frac{L_f^2}{3\nu},
$$
the kinetic energy grows at most linearly in time
$$
E(t) \le 2E(t_0) + F^2\left (\frac{(t-t_0)L_f^2}{\nu}\right ).
$$
The long time for the average dissipation rate is bounded
\begin{equation}
\lim\sup_{t\to\infty} \left < \epsilon(\cdot )\right >_t \le \frac{F^2L_f^2}{\nu} = \epsilon_B. \label{epsilonb}
\end{equation}
These bounds are uniform in the size $L$ of the period which we assume to be
much larger than $L_f$. If the size of the period is allowed to enter the calculations then the kinetic energy is bounded by
$$
E(t) \le L^2 \frac{L_f^2F_*^2}{\nu^2} + \left (E(t_0) -
L^2\frac{L_f^2F_*^2}{\nu^2}\right )e^{-\frac{\nu(t-t_0)}{L^2}}
$$
where 
$$
L_f^2F_*^2 = \sup_t L^{-3}\int |(-\Delta)^{-\frac{1}{2}}f(x,t)|^2dx.
$$
This means that for much longer times
$$
t-t_0\ge \frac{L^2}{\nu}
$$
the kinetic energy saturates to a value that depends on the large scale. 
But the bound (\ref{epsbound}) that is independent of $L$ is always valid; it can be written in
terms of
\begin{equation}
B = 4E(t_0) + (t-t_0)\epsilon_B
\label{energy}
\end{equation}
as
\begin{equation}
E(t) + (t-t_0)\left <\epsilon \right >_t \le  B. 
\label{energyb}
\end{equation}
A useful K-bound is
\begin{equation}
\int_{t_0}^t\|u(\cdot,s)\|_{L^{\infty}(dx)}ds \le K_{\infty}\label{maxu}
\end{equation}
The constant $K_{\infty}$ has dimensions of length and depends
on the initial kinetic energy, viscosity, body forces and time. The bound follows
by interpolation from (\cite{gft}) and is derived in Appendix A
together with the formula 
\begin{equation}
K_{\infty} = C\left\{ \frac{K_0}{\nu^2} + \sqrt{\nu (t-t_0)} +
\frac{t-t_0}{\nu^2}\int_{t_0}^t\|f(\cdot,s)\|_{L^2}^2ds\right\}.
\label{kif}
\end{equation}
The displacement $\ell$ satisfies certain K-bounds that follow from
the bounds above and (\ref{deltaeq}). We mention here
\begin{equation}
\|\ell(\cdot ,t)\|_{L^{\infty}(dx)}
\le \int_{t_0}^t\|u(\cdot,s)\|_{L^{\infty}(dx)}ds \le K_{\infty},\label{maxdel}
\end{equation}
The inequality (\ref{maxdel}) follows from (\ref{deltaeq}) by multiplying with $\ell |\ell |^{2(m-1)}$,
integrating,
$$
\frac{1}{2m}\frac{d}{dt}\int |\ell(x,t)|^{2m}dx + \nu\int|\nabla\ell(x,t)|^2|\ell(x,t)|^{2(m-1)}dx +
$$
$$
+ \nu\frac{m-1}{2}\int |\nabla|\ell (x,t)|^2|^2|\ell (x,t)|^{2(m-2)}dx  +
$$
\begin{equation}
+ \int u(x,t)\cdot\ell(x,t) |\ell(x,t)|^{2(m-1)} dx \le 0,\label{deltum}
\end{equation}
and then ignoring the viscous terms, using H\"{o}lder's inequality
in the last term, multiplying by $m$, 
taking the $m$-th root, integrating in time and then letting $m\to \infty$. 

The case $m=1$ gives
$$
\frac{d}{2dt}\int |\ell(x,t)|^2dx + \nu \int|\nabla \ell(x,t)|^2 dx
\le \sqrt{K_0}\sqrt{\int |\ell(x,t)|^2dx}
$$
and consequently, we obtain by integration  from $t_0 =0$
\begin{equation}
\sqrt{\int|\ell (x,t)|^2dx} \le  
t\sqrt{K_0}, \label{elltwo}
\end{equation}
and then, using (\ref{elltwo}) we deduce the inequality
\begin{equation}
\int_0^t\int |\nabla \ell(x,s)|^2 dxds \le \frac{{K_0}t^2}{2\nu}.
\label{nablaelltwo}
\end{equation}
Now we  multiply (\ref{deltaeq}) by $-\Delta \ell$, 
integrate by parts, use Schwartz's inequality to write
$$
\frac{d}{dt}\int |\nabla\ell (x,t)|^2dx + 
\nu\int |\Delta\ell (x,t)|^2dx \le
\sqrt{\int|\nabla u(x,t)|^2dx}{\sqrt{\int |\nabla \ell(x,t)|^2dx}}
$$
$$
  - 2\int{\mbox Trace}\left \{(\nabla\ell(x,t))(\nabla u(x,t))(\nabla \ell(x,t))^*\right \}dx
$$
and then use the elementary inequality
$$
\left (\int |\nabla\ell (x,t)|^4dx\right )^{\frac{1}{2}} \le C\|\ell (\cdot,t)\|_{L^{\infty}}
\left (\int|\Delta \ell (x,t)|^2dx\right )^{\frac{1}{2}},
$$
in conjunction with the H\"{o}lder inequality and (\ref{maxdel}) to deduce
$$
\frac{d}{dt}\int |\nabla\ell (x,t)|^2dx + 
\nu\int |\Delta\ell (x,t)|^2dx \le
$$
$$
\sqrt{\int|\nabla \ell (x,t)|^2dx}\sqrt{\int |\nabla u (x,t)|^2dx} +
C\frac{K_{\infty}^2}{\nu}\int |\nabla u (x,t)|^2dx.
$$
We obtain, after integration and use of (\ref{en}, \ref{nablaelltwo})

\begin{equation}
\int |\nabla \ell(x,t)|^2 dx + \nu \int_{0}^t\int |\Delta\ell(x,s)|^2dxds \le C\left (\frac{K_0t}{\nu} +
\frac{K_{\infty}^2K_0}{\nu^2}\right ).
\label{deldel}
\end{equation} 
Recalling the bound (\ref{energy}, \ref{energyb}) on kinetic energy we have:
\begin{thm} Assume that the vector valued function $\ell$ obeys
(\ref{deltaeq}) and assume that the velocity $u(x,t)$ is a
solution of the Navier-Stokes equations (or, more generally, that it 
is a divergence-free periodic function that satisfies
the bounds (\ref{en}) and (\ref{maxu})). Then $\ell$ satisfies
the inequality (\ref{maxdel}) together with
\begin{equation}
\frac{1}{L^3}\int|\ell (x,t)|^2 dx \le  
(4E_0  + t\epsilon_B) t^2, 
\label{ltwo}
\end{equation}
\begin{equation}
\frac{1}{L^3 t}\int_0^t\int |\nabla \ell(x,s)|^2 dx ds \le 
\frac{Bt}{2\nu},
\label{nablaeltwo}
\end{equation}
and
\begin{equation}
\int |\nabla \ell(x,t)|^2 \frac{dx}{L^3} + \nu \int_{0}^t\int |\Delta\ell(x,s)|^2\frac{dx}{L^3}ds\le C\left (\frac{B t}{\nu} +
\frac{K_{\infty}^2B}{\nu^2}\right ).
\label{deltaltwo}
\end{equation}
In these inequalities 
$$
E_0 = \frac{1}{2L^3}\int |u(x,0)|^2dx,
$$
$$
B = 4E_0 + t\epsilon_B
$$
and $\epsilon_B$ is given in (\ref{epsilonb}).
\end{thm}

Let us consider the pair dispersion
\begin{equation}
\left <\delta_t^2\right > = L^{-6}\int\int_{\{(x,y); |A(x,t)-A(y,t)|\le \delta_0\}}
|x-y|^2 dxdy.\label{pair}
\end{equation}
Using the triangle inequality (\ref{tri}) in (\ref{ltwo}) we obtain
\begin{thm}
Consider periodic solutions of the Navier-Stokes equation
with large period $L$, and assume that the body forces have $L_f$
finite. Then the
pair dispersion obeys
\begin{equation}
\left <\delta_t^2\right> \le 3\delta_0^2 + 24t{E_0}t^2  + 
6\epsilon_B t^3.\label{richa}
\end{equation}
\end{thm}

\noindent {\bf Comment} Use
of the ODE $\frac{dX}{dt} = u(X,t)$ requires information about the
gradient $\nabla A$ and produces worse bounds.

\section{$\epsilon$-bounds}
This section is devoted to bounds on higher order derivatives
of $\ell$. These bounds require assumptions.
We are going to apply the Laplacian to
(\ref{deltaeq}), multiply by $\Delta\ell $ and integrate. We obtain
$$
\frac{1}{2}\frac{d}{dt}\int|\Delta \ell (x,t)|^2dx + \nu\int |\nabla\Delta\ell (x,t)|^2dx =
$$
\begin{equation}
\int \partial_k u(x,t)\cdot\partial_k\Delta\ell (x,t)dx +
I
\label{inte}
\end{equation}
where
$$
I = \int \partial_k(u(x,t)\cdot\nabla \ell(x,t))\cdot\partial_k\Delta\ell(x,t)dx.
$$
Now 
$$
I = \int(\partial_k u)\cdot\nabla\ell (x,t)\cdot \partial_k\Delta\ell (x,t)dx + II
$$
where
$$
II = \int u(x,t)\cdot\nabla(\partial_k\ell(x,t))\cdot\Delta\partial_k\ell (x,t)dx
$$
and, integrating by parts
$$
II = - \int\partial_lu(x,t)\cdot\nabla (\partial_k\ell(x,t))\cdot\partial_l\partial_k\ell(x,t)dx
$$
and then again
$$
II = \int\partial_lu(x,t)\cdot\nabla \partial_l\partial_k\ell(x,t) 
\cdot(\partial_k\ell (x,t))dx
$$
and so 
$$
I = \int\partial_l u_i(x,t)\partial_k \ell_j(x,t)\left \{\partial_i\partial_k
 + \delta_{ik}\Delta \right \} \partial_l\ell_j(x,t)dx
$$
Putting things together we get
$$
|I|\le C\|\nabla\ell(\cdot, t)\|_{L^{\infty}}\|\nabla u(\cdot,t)\|_{L^2}\|\nabla\Delta\ell(\cdot, t)\|_{L^2}
$$
Thus
$$
\frac{1}{2}\frac{d}{dt}\int|\Delta \ell (x,t)|^2dx + \nu\int |\nabla\Delta\ell
(x,t)|^2dx \le
$$
\begin{equation}
\frac{C}{\nu}\int|\nabla u(x,t)|^2dx +  C\|\nabla\ell(\cdot, t)\|_{L^{\infty}}\|\nabla u(\cdot,t)\|_{L^2}\|\nabla\Delta\ell(\cdot, t)\|_{L^2}\label{ine}
\end{equation}
Now we use an interpolation inequality that is valid for periodic functions
with zero mean and implies that
$$
\|\nabla \ell\|_{L^{\infty}} \le c \|\Delta \ell\|_{L^2}^{\frac{1}{2}}\|\nabla\Delta\ell\|_{L^2}^{\frac{1}{2}}
$$
Using this inequality we obtain
$$
\frac{d}{dt}\int|\Delta \ell (x,t)|^2dx + \nu\int |\nabla\Delta\ell
(x,t)|^2dx \le
$$
\begin{equation}
\frac{C}{\nu}\int|\nabla u(x,t)|^2dx +  C\nu^{-3}\|\nabla u(\cdot,t)\|_{L^2}^4\|\Delta\ell(\cdot, t)\|_{L^2}^2\label{ineq}
\end{equation}
Therefore we deduce
\begin{equation}
\int |\Delta \ell(x,t)|^2\frac{dx}{L^3}\le c\frac{B}{\nu^2}{\mbox{exp}}\left
\{{\frac{c L^6}{\nu^5}\int_0^t\epsilon^2(s)ds}\right\}
\label{ug}
\end{equation}
where
\begin{equation}
\epsilon (s) = \nu L^{-3}\int|\nabla u(x,s)|^2dx
\label{eps}
\end{equation}
is the instantaneous energy dissipation.
\begin{prop} If $\ell$ solves (\ref{deltaeq}) with periodic boundary 
conditions on a time interval $t\in [0,T]$ and if the integral
$$
\int_0^T\epsilon^2 (s)ds
$$
is finite, then
$$
\int  |\Delta \ell(x,t)|^2\frac{dx}{L^3} + \nu \int_0^t\int |\nabla \Delta
\ell (x,s)|^2\frac{dx}{L^3} \le c \frac{B}{\nu^2}{\mbox{exp}}\left
\{{\frac{c L^6}{\nu^5}\int_0^t\epsilon^2(s)ds}\right\}
$$
holds for all $0\le t\le T$.
\end{prop}

\section{Bounds for the virtual velocity}
We prove here the assertion that $v$ does not grow too much as
long as the $L^3$ norm of $C$ is not too large.
We recall that $v$ solves (\ref{veq})
\begin{equation}
\Gamma v_i = 2\nu C_{m,k;i}\partial_kv_m + Q_{ji}f_j.\label{ve}
\end{equation}
We multiply by $v_i |v|^{2(m-1)}$ and integrate:
$$
\frac{1}{2m}\frac{d}{dt}\int |v(x,t)|^{2m}dx + \nu\int|\nabla v(x,t)|^2|v(x,t)|^{2(m-1)}dx +
$$
$$
+ \nu\frac{m-1}{2}\int |\nabla|v(x,t)|^2|^2|v(x,t)|^{2(m-2)}dx   =
$$
$$
= 2\nu\int C_{m,k;i}(x,t)(\partial_kv_m(x,t))v_i(x,t) |v(x,t)|^{2(m-1)}dx +
$$
\begin{equation}
+ \int  Q_{ji}(x,t)f_j(x,t)v_i(x,t) |v(x,t)|^{2(m-1)}dx.\label{vum}
\end{equation}

We bound
$$
2\nu\left |\int C_{m,k;i}(x,t)(\partial_kv_m(x,t))v_i(x,t) |v(x,t)|^{2(m-1)}dx \right | \le 
$$
$$
\nu\int|\nabla v(x,t)|^2|v(x,t)|^{2(m-1)}dx + \nu \int|C(x,t)|^2|v(x,t)|^{2m}dx
$$
where
\begin{equation}
|C(x,t)|^2 = \sum_{m,k,i} |C_{m,k;i}(x,t)|^2, \label{C}
\end{equation}
and we bound
$$
\left |\int  Q_{ji}(x,t)f_j(x,t)v_i(x,t) |v(x,t)|^{2(m-1)}dx\right |\le
$$
$$
\left \{\int |g(x,t)|^{2m}dx\right \}^{\frac{1}{2m}}
\left \{\int |v(x,t)|^{2m}dx \right \}^{\frac{2m-1}{2m}}
$$
where
\begin{equation}
g_i(x,t) = Q_{ji}(x,t)f_j(x,t).\label{g}
\end{equation}
The inequality obtained is
$$
\frac{d}{dt}\int |v(x,t)|^{2m}dx + \nu m(m-1)\int |\nabla|v(x,t)|^2|^2|v(x,t)|^{2(m-2)}dx \le
$$
$$
\le 2m\nu \int|C(x,t)|^2|v(x,t)|^{2m}dx\, + 
$$
\begin{equation}
+ \, 2m\left \{\int |g(x,t)|^{2m}dx\right \}^{\frac{1}{2m}}
\left \{\int |v(x,t)|^{2m}dx \right \}^{\frac{2m-1}{2m}}\label{vineq}
\end{equation}
Let us consider for any $m\ge 1$ the quantity 
$$
q(x,t) = |v(x,t)|^m.
$$
The inequality (\ref{vineq}) implies
$$
\frac{d}{dt}\int (q(x,t))^2 + 4\nu(1-\frac{1}{m})\int |\nabla q (x,t)|^2
\le 2m\nu \int|C(x,t)|^2(q(x,t))^2dx\, + 
$$
$$
+ \, 2m\left \{\int |g(x,t)|^{2m}dx\right \}^{\frac{1}{2m}}
\left \{\int (q(x,t))^2dx \right \}^{\frac{2m-1}{2m}}
$$
Using the well-known Morrey-Sobolev inequality
$$
\left \{\int (q(x))^6dx\right \}^{\frac{1}{3}} \le C_0\left \{ \int |\nabla q(x)|^2dx + L^{-2}\int (q(x))^2dx\right \}
$$
and H\"{o}lder's inequality we deduce
$$
\frac{d}{dt}\int (q(x,t))^2 + 4\nu(1-\frac{1}{m})\int |\nabla q (x,t)|^2
\le
$$
$$
\le 2m\nu C_0\left \{\int |C(x,t)|^3dx\right\}^{\frac{2}{3}} \left \{ \int |\nabla q(x,t)|^2dx + L^{-2}\int (q(x,t))^2dx\right \}
$$

$$
+ \, 2m\left \{\int |g(x,t)|^{2m}dx\right \}^{\frac{1}{2m}}
\left \{\int (q(x,t))^2dx \right \}^{\frac{2m-1}{2m}}
$$
Therefore, if on the time interval $t\in [0, \tau]$ $C(x,t)$ obeys the 
smallness condition
\begin{equation}
\left \{\int |C(x,t)|^3dx\right\}^{\frac{1}{3}} \le \sqrt{\frac{2(m-1)}{C_0m^2}}\label{Ccond}
\end{equation}
then we have the inequality
$$
\frac{d}{dt}\|v(\cdot, t)\|_{L^{2m}}\le \frac{\nu(m-1)}{2m^2L^2}
\|v (\cdot, t)\|_{L^{2m}} + \|g(\cdot, t)\|_{L^{2m}}
$$
for $t\in [0,\tau]$ and consequently
\begin{equation}
\|v(\cdot ,t)\|_{L^{2m}} \le \|v_0\|_{L^{2m}}e^{\frac{\nu(m-1)t}{2m^2L^2}}
+ \int_0^t\|g(\cdot, s)\|_{L^{2m}} \label{vbound}
\end{equation}
holds on the same time interval.

\section{Appendix A}
In this appendix we prove the inequality (\ref{maxu}) and derive
the explicit expression for $K_{\infty}$. The calculation is based 
on (\cite{gft}). All constants $C$ are non-dimensional and
may change from line to line. Solutions $u$ of the Navier-Stokes equations 
obey the differential inequality
$$
\frac{d}{ds}\int|\nabla u(x,s)|^2dx + \nu\int|\Delta u (x,s)|^2dx
$$
$$
\le \frac{C}{\nu^3}\left (\int |\nabla u(x,s)|^2dx\right )^3 +
\frac{C}{\nu}\int |f(x,s)|^2dx.
$$
The idea of (\cite {gft}) was to divide by an appropriate quantity to make use of the balance (\ref{en}). The quantity is
$$
(G(s))^2 = \left (\gamma^2 + \int|\nabla u(x,s)|^2dx\right )^{2}
$$
where $\gamma$ is a positive constant that does not depend
on $s$ and will be specified later. Dividing by $(G(s))^2$,
integrating in time from $t_0$ to $t$ and using (\ref{en}) one obtains
$$
\int_{t_0}^t \|\Delta u (\cdot, s)\|^2_{L^2}(G(s))^{-2}ds \le 
$$
$$
C\left(\frac{K_0}{\nu^5} + \frac{1}{\nu\gamma^2} +
\frac{1}{\nu^2\gamma^4}\int_{t_0}^t\|f(\cdot,s)\|^2ds\right ).
$$
The three dimensional Sobolev embedding-interpolation inequality 
for periodic mean-zero functions
$$
\|u\|_{L^{\infty}} \le C\|\nabla u\|_{L^2}^{\frac{1}{2}}\|\Delta u\|_{L^2}^{\frac{1}{2}}
$$ 
is elementary. From it we deduce
$$
\|u(\cdot,s)\|_{L^{\infty}} \le C \|\nabla u(\cdot, s)\|_{L^2}^{\frac{1}{2}}(G(s))^{\frac{1}{2}}\left [\|\Delta u(\cdot,s)\|_{L^2}G(s)^{-1}\right ]^{\frac{1}{2}}
$$
Integrating in time, using the H\"{o}lder inequality, the inequality 
(\ref{en}) and the inequalities above we deduce
$$
\int_{t_0}^t\|u(\cdot,s)\|_{L^{\infty}} \le Cr 
$$
where the length $r=r(\gamma,t,\nu,K_0)$ is given in terms of six length scales
$$
\frac{K_0}{\nu^2} = r_0, \quad \frac{\nu^2}{\gamma^2} = r_1, \quad
\frac{(t-t_0)\gamma^2}{\nu} = r_2,
$$
$$
(\gamma (t-t_0))^{\frac{2}{3}} = r_3, \quad \frac{t-t_0}{\nu^2}\int_{t_0}^t\|f(\cdot,s)\|_{L^2}^2ds = r_4
$$
and
$$
r_5 =\sqrt{\nu (t-t_0)}. 
$$
The expression for $r$ is
$$
r = r_0 + (r_0)^{\frac{3}{4}} (r_1)^{\frac{1}{4}} +  (r_0)^{\frac{1}{2}} (r_2)^{\frac{1}{2}} + (r_0)^{\frac{1}{4}} (r_3)^{\frac{3}{4}} +
$$
$$
(r_0)^{\frac{1}{4}} (r_4)^{\frac{1}{4}}(r_5)^{\frac{1}{2}} + (r_0)^{\frac{3}{4}} (r_4)^{\frac{1}{4}}\left (\frac{r_1}{r_2}\right)^{\frac{1}{4}}
$$
The choice
$$
\gamma ^4 = \frac{\nu^3}{t-t_0}
$$
entrains
$$
r_1 = r_2 = r_3 = r_5
$$
reducing thus the number of length scales to three, the energy viscous
length scale $r_0$, the diffusive length scale $r_5$ and the force
length scale $r_4$. The bound becomes 
$$
K_{\infty} = C(r_0 + r_4 + r_5)
$$
i.e. (\ref{kif}).

\section{Appendix B}
We prove her the commutation relation (\ref{lagcom}). We take an arbitrary
function $g$ and compute $[\Gamma, L_i g])$ where 
$\Gamma = \Gamma_{\nu}(u,\nabla)$ and $L_i = \nabla_A^i$. We use first (\ref{lfg}):
$$
[\Gamma, L_i g] = \Gamma\left (Q_{ji}\partial_j g\right ) - Q_{ji}\partial_j\Gamma g =
$$ 
$$
\Gamma(Q_{ji}) \partial_j g  + Q_{ji}\Gamma\partial_j g -
2\nu\partial_k(Q_{ji})\partial_k\partial_j g - Q_{ji}\partial_j\Gamma g =
$$
(commuting in the last term $\partial_j$ and $\Gamma$)
$$
\Gamma(Q_{ji}) \partial_j g - 2\nu\partial_k(Q_{ji})\partial_k\partial_j g
-Q_{ji}\partial_j(u_k)\partial_k g =
$$
(changing names of dummy indices in the last term)
$$
\left (\Gamma (Q_{ji}) - Q_{ki}\partial_k (u_j)\right )\partial_j g
-  2\nu\partial_k(Q_{ji})\partial_k\partial_j g =
$$
(using (\ref{Qeq}))
$$
2\nu Q_{jp}(\partial_l\partial_k A_p)(\partial _kQ_{li})\partial_j g -
2\nu\partial_k(Q_{ji})\partial_k\partial_j g =
$$
(using the definition (\ref{lai}) of $\nabla_A$
$$
2\nu(\partial_l\partial_k A_p)(\partial _kQ_{li})(L_p g) - 
2\nu\partial_k(Q_{ji})\partial_k\partial_j g =
$$
(renaming dummy indices in the last expression)
$$
2\nu (\partial_k(Q_{li}))\left \{(\partial_l\partial_k A_p)(L_p g) -
\partial_l\partial_k g\right \} = 
$$
(using (\ref{ial}) in the last expression)
$$
2\nu (\partial_k(Q_{li}))\left \{(\partial_l\partial_k A_p)(L_p g) -
\partial_k \left (\partial_l(A_p)(L_p g)\right )\right \} =
$$
(carrying out the differentiation in the last term and cancelling)
$$
-2\nu(\partial_k(Q_{li}))\partial_l(A_p)\partial_k(L_p(g)) =
$$ 
(using the differential consequence of the
fact that $Q$ and $\nabla A$ are inverses of each other)
$$
2\nu Q_{li}(\partial_k\partial_l(A_p)) \partial_k(L_p g) =
$$
(using the definition (\ref{lai}) of $\nabla_A$) 
$$
2\nu L_i(\partial_k(A_p))\partial_k(L_p g) =
$$
(using the definition (\ref{cmki}) of $C_{m,k;i}$)
$$
2\nu C_{p,k; i}\partial_k(L_p g),
$$
and that concludes the proof. 
We proceed now to prove (\ref{curveq}). We start with  (\ref{nablaeq})
$$
\Gamma(\partial_k A_m) = - (\partial_k u_j)(\partial_j A_m)
$$
and apply $L_i$:
$$
L_i(\Gamma(\partial_k A_m)) = - L_i\left \{(\partial_k u_j)(\partial_j A_m)\right \}.
$$
Using the commutation relation (\ref{lagcom}) and the definition (\ref{cmki})
we get
$$
\Gamma (C_{m,k;i}) = - L_i\left \{(\partial_k u_j)(\partial_j A_m)\right \} +
2\nu C_{p,l;i}\partial_lL_p(\partial_k A_m).
$$
Using the fact that $L_i$ is a derivation in the first term
and the definition  (\ref{cmki}) in the last term we conclude that
$$
\Gamma (C_{m,k;i}) = -(\partial_k u_j)C_{m,j;i} - (\partial_j A_m)(L_i(\partial_k u_j)) + 2\nu C_{p,l;i}(\partial_l C_{m,k; p})
$$ 
which is (\ref{curveq}). We compute now the formal adjoint of $\nabla_A^i$
$$
(\nabla_A^i)^*g = -\partial_j(Q_{ji}g) = -\nabla_A^i(g) - (\partial_j(Q_{ji}))g
$$
(with (\ref{ial}))
$$
(\nabla_A^i)^*g = - \nabla_A^i(g) - \left \{(\partial_jA_p)L_p(Q_{ji})\right \}g =
$$
(using the fact that $Q$ is the inverse of $\nabla A$)
$$
(\nabla_A^i)^*g = - \nabla_A^i(g) + Q_{ji}C_{p,j;p}g.
$$

{\bf Acknowledgments.} Part of this work was done at
the Institute for Theoretical Physics
in Santa Barbara, whose hospitality is gratefully acknowledged.
This research is supported in part by NSF- DMS9802611.


\begin{thebibliography}{99}
\bibitem{cel} P. Constantin, An Eulerian-Lagrangian approach for incompressible fluids: Local theory, http://arXiv.org/abs/math.AP/0004059.
\bibitem{celw} P.Constantin, An Eulerian-Lagrangian approach to fluids, preprint 1999, (www.aimath.org). 
\bibitem{dc} C. Doering, P. Constantin,  Energy dissipation in shear driven turbulence, Phys.Rev.Lett. {\bf 69} (1992), 1648-1651.
\bibitem{moya} A. S. Monin and A. M. Yaglom, {\em Statistical
 Fluid Mechanics}
M.I.T. Press,  Cambridge, 1987.
\bibitem{frisch} U. Frisch,  {\em Turbulence}, Cambridge Univ
ersity Press, Cambridge,
1995.
\bibitem{labexp} J. Mann, S. Ott and J.S. Andersen, Experimental study
of relative turbulent diffusion, Riso National Laboratory, Denmark
Riso-R-1036(EN) (1999).
\bibitem{serr} J. Serrin, Mathematical principles of classica
l fluid
mechanics, (S. Flugge,
C. Truesdell Edtrs.) Handbuch der Physik, {\bf 8} (1959), 125
-263,
p.169.
\bibitem{gold} M. E. Goldstein, Unsteady vortical and entropic
distortion of potential flows round arbitrary obstacles, J. Fluid
Mech. {\bf 89} (1978), 433-468.
\bibitem{goldet} M. E. Goldstein, P. A. Durbin, The effect of
 finite turbulence spatial scale on the amplification of 
turbulence by a contracting stream, J. Fluid Mech. {\bf 98} (1980), 473-508. 
\bibitem{hunt}  J. C. R. Hunt, Vorticity and vortex dynamics 
in complex turbulent flows, Transactions of CSME, {\bf 11} (1987), 21-35.
\bibitem{kuz} G. A. Kuzmin, Ideal incompressible hydrodynamics
in terms of the vortex momentum density, Phys. Lett {\bf 96 A
} (1983), 88-90.
\bibitem{ose} V. I. Oseledets, On a new way of writing the
  Navier-Stokes equation. The Hamiltonian formalism, Commun. 
Moscow Math. Soc. (1988), Russ. Math. Surveys {\bf 44} (1989), 210
-211.
\bibitem{h} Shiyi Chen, Ciprian Foias, Darryl D Holm, Eric Olson, 
Edriss S Titi, Shannon Wynn, A connection between the
Camassa-Holm equations and 
turbulent flows in channels and pipes, Phys. Fluids, {\bf{11}}
(1998), 2343-2353.
\bibitem{mar} D.D. Holm, J. E. Marsden, T. Ratiu, Euler-Poincar\'{e} 
models of ideal fluids with nonlinear dispersion, Phys. Rev. Lett.
{\bf 349} (1998), 4173-4177.
\bibitem{cfbook} P. Constantin, C. Foias, Navier-Stokes equations,
University of Chicago Press, Chicago, 1988.
\bibitem{gft} C. Guillop\'{e}, C. Foias, R. Temam, New a priori estimates for the Navier-Stokes equations in dimension 3, Commun. PDE {\bf 6} (1981),329-359.
\end{thebibliography}
\end{document}